# Dispersive Divergence-Free Vector Meshless Method for Analysis of Frequency-Dependent Medium

Sheyda Shams
*Department of Electrical and Computer Engineering, Ph.D., Candidate*
Shiraz University
Shiraz, Iran
sheyda.shams@shirazu.ac.ir

Masoud Movahhedi
*Department of Electrical and Computer Engineering, Ph.D., Associate Professor*
Yazd University
Yazd, Iran
movahhedi@yazd.ac.ir

*Abstract*—The dispersive meshless method with scalar basis function has been successfully applied for analysis of frequency-dependent media. However, as scalar-based meshless methods are not always divergence-free in the absence of source, inaccurate or even wrong solutions may be found in the results. To overcome this problem, in this paper, a new dispersive formulation of meshless method with vector basis function is proposed for analysis of dispersive materials. The divergence-free property of this method improves the accuracy of simulation by eliminating the artificial charges generated due to the numerical analysis. Moreover, the frequency behavior of the dispersive medium is modeled using auxiliary-differential-equations (ADE) method. The efficiency and divergence-free property of the proposed method are verified and compared with dispersive meshless method with scalar basis function by a numerical example.

## I. Introduction

Numerical analysis of dispersive mediums is of considerable importance because the constitutive parameters of many real materials are dependent on frequency. Conventional transient numerical methods are efficient and powerful tools for modeling of electromagnetic waves propagation in mediums which their constitutive parameters are not frequency-dependent. However, when it comes to numerical analysis of dispersive mediums, the conventional numerical methods suffer from inability in modeling of frequency-dependent properties of dispersive media.

Since 1990, a number of dispersive formulation for grid-based numerical methods such as finite-difference time-domain (FDTD) and finite-element method (FEM) are proposed. Some of these works can be found in [1]-[4]. Although FDTD and FEM are well-developed numerical techniques, they do not possess the capability of meshless methods in modeling complex structures. In meshless methods, spatial domain of problems is discretized by a set of nodes without obtaining connection information among nodes. This considerable property in addition to their high accuracy make meshless methods attractive for analysis of dispersive media [5].

In [5], for the first time, a dispersive meshless formulation for analyzing frequency-dependent materials was presented. That formulation was based on the scalar radial basis function (RBF) meshless method [6] for discretizing the scalar form of the Maxwell equations in spatial domain and application of auxiliary differential equation (ADE) method for modeling the frequency behavior of dispersive medium. In [5], the effectiveness of the proposed dispersive meshless method in comparison to the grid-based numerical methods for analyzing dispersive mediums was investigated and verified by numerical examples.

Based on Maxwell's equations, the divergence of magnetic fields are always equal to zero and electric fields are divergence-free in source-free regions. However, due to discretization scheme in numerical analysis of electromagnetic problems, the divergence-free property is not always preserved, and the accuracy of simulations may be reduced. To address this issue, [7] presented a divergence-free vector-based meshless method (vector meshless method) for solving Maxwell's equation. References [8]-[10] include number of studies in which this method was effectively applied for analyzing propagation of electromagnetic waves.

Up to now based on our knowledge, no reports on applications of the vector meshless method in analysis of dispersive mediums were found. In this paper, to obtain a divergence-free formulation for dispersive meshless method, we have used meshless method with vector basis functions for discretizing the vector form of Maxwell's equations. Moreover, ADE technique is applied to vector meshless method to include the frequency behavior of a dispersive medium, such as plasma or metamaterial, with Drude dispersion model [11]. We show that the divergence-free property of the proposed method leads to eliminate the spurious charges in source-free regions of problem domain and improves the accuracy of simulations.

## II. Dispersive Formulation of Divergence-free Vector-based Meshless method

For modeling the frequency behavior of the dispersive media, we apply ADE method in this paper because its implementation is simpler than the other existing approaches according to [3]. Among different algorithms of ADE method, we use auxiliary differential equations based on flux-density, which relates the electric and magnetic flux-density vectors, *i.e.*, $\mathbf{D}$ and $\mathbf{B}$, to the electric and magnetic field-intensity vectors, *i.e.*, $\mathbf{E}$ and $\mathbf{H}$, respectively. To obtain the dispersive formulation of vector meshless method, the vector form of Maxwell's equations with the vector of electric flux-density $\mathbf{D}$ and magnetic flux-density $\mathbf{B}$ is considered as follows:

$$\frac{\partial \mathbf{D}}{\partial t} = \nabla \times \mathbf{H} \qquad (1)$$

$$\frac{\partial \mathbf{B}}{\partial t} = -\nabla \times \mathbf{E} \qquad (2)$$

We discretize the spatial domain using a set of **D** and a set of **B** nodes in addition to the **E** and **H** nodes [5]. Based on the definition of vector shape function $\Phi_j$ [7] in vector meshless method, the electromagnetic fields and flux densities can be approximated using following equations:

$$\mathbf{H} = \sum_j^{NH} \Phi_j \mathbf{H}_j \quad (3)$$

$$\mathbf{D} = \sum_j^{ND} \Phi_j \mathbf{D}_j \quad (4)$$

$$\mathbf{B} = \sum_j^{NB} \Phi_j \mathbf{B}_j \quad (5)$$

$$\mathbf{E} = \sum_j^{NE} \Phi_j \mathbf{E}_j \quad (6)$$

By substitution of (3)-(6) into (1) and (2), we obtain the following equations:

$$\frac{\partial \sum_j^{ND} \Phi_j \mathbf{D}_j}{\partial t} = \nabla \times \sum_j^{NH} \Phi_j \mathbf{H}_j \quad (7)$$

$$\frac{\partial \sum_j^{NB} \Phi_j \mathbf{B}_j}{\partial t} = -\nabla \times \sum_j^{NE} \Phi_j \mathbf{E}_j \quad (8)$$

Acording to Kronecker's delta property of vector shape function [7], and using central difference scheme for approximating the time-derivatives of Maxwell's equations, the discretized form of Maxwell's equations is obtained as follows:

$$\mathbf{B}^{n+3/2} = \mathbf{B}^{n+1/2} - \Delta t \sum_j^{NE} \nabla \times \Phi_j \mathbf{E}_j^{n+1} \quad (9)$$

$$\mathbf{D}^{n+1} = \mathbf{D}^n + \Delta t \sum_j^{NH} \nabla \times \Phi_j \mathbf{H}_j^{n+1/2} \quad (10)$$

Based on Drude dispersion model, the electromagnetic flux densities are related to electromagnetic field intensities as follows [5]:

$$\mathbf{D} = \varepsilon_0 (1 - \frac{\omega_{ep}^2}{\omega^2 - j\gamma_e \omega}) \mathbf{E} \quad (11)$$

$$\mathbf{B} = \mu_0 (1 - \frac{\omega_{mp}^2}{\omega^2 - j\gamma_m \omega}) \mathbf{H} \quad (12)$$

where $\omega_{ep}$ ($\omega_{mp}$) is the electric (magnetic) plasma frequency and $\gamma_e$ ($\gamma_m$) is the corresponding collision frequency [12].

As the inverse Fourier transform of $j\omega$ is equal to the time-derivatives in the time-domain, obtaining the inverse Fourier transform of dispersion equations of (11) and (12) leads to the following equations in the time-domain:

$$\frac{\partial^2 \mathbf{D}}{\partial t^2} + \gamma_e \frac{\partial \mathbf{D}}{\partial t} = \varepsilon_0 \frac{\partial^2 \mathbf{E}}{\partial t^2} + \varepsilon_0 \gamma_e \frac{\partial \mathbf{E}}{\partial t} + \varepsilon_0 \omega_{ep}^2 \mathbf{E} \quad (13)$$

$$\frac{\partial^2 \mathbf{B}}{\partial t^2} + \gamma_m \frac{\partial \mathbf{B}}{\partial t} = \mu_0 \frac{\partial^2 \mathbf{H}}{\partial t^2} + \mu_0 \gamma_m \frac{\partial \mathbf{H}}{\partial t} + \mu_0 \omega_{mp}^2 \mathbf{H} \quad (14)$$

We approximate the first and second order time derivatives of **E** and the other vector quantities in (13) and (14) using following central-difference scheme centered at time-step *n* as [12]:

$$\frac{\partial \mathbf{E}}{\partial t} = \frac{\mathbf{E}^{n+1} - \mathbf{E}^{n-1}}{2\Delta t} \quad (15)$$

$$\frac{\partial^2 \mathbf{E}}{\partial t^2} = \frac{\mathbf{E}^{n+1} - 2\mathbf{E}^n + \mathbf{E}^{n-1}}{(\Delta t)^2} \quad (16)$$

Moreover, **E** and **H** fields are approximated by a semi-implicit scheme as [12]:

$$\mathbf{E} = \frac{\mathbf{E}^{n+1} + 2\mathbf{E}^n + \mathbf{E}^{n-1}}{4} \quad (17)$$

Finally, the dispersion equations of (11) and (12) in the time-domain are formulated as:

$$\mathbf{H}^{n+3/2} = \sum_{m=0}^{M} d_m \mathbf{B}^{n-m+3/2} - \sum_{m=1}^{M} c_m \mathbf{H}^{n-m+3/2} \quad (18)$$

$$\mathbf{E}^{n+1} = \sum_{m=0}^{M} b_m \mathbf{D}^{n-m+1} - \sum_{m=1}^{M} a_m \mathbf{E}^{n-m+1} \quad (19)$$

where $M = 2$ and coefficients are expressed as:

$$a_1 = \frac{2\varepsilon_0 \Delta t^2 \omega_{ep}^2 - 8\varepsilon_0}{A} \quad (20)$$

$$a_2 = \frac{4\varepsilon_0 - 2\Delta t \varepsilon_0 \gamma_e + \varepsilon_0 \Delta t^2 \omega_{ep}^2}{A} \quad (21)$$

$$b_0 = \frac{4 + 2\Delta t \gamma_e}{A} \quad (22)$$

$$b_1 = \frac{-8}{A} \quad (23)$$

$$b_2 = \frac{4 - 2\Delta t \gamma_e}{A} \quad (24)$$

$$A = 4\varepsilon_0 + 2\Delta t \varepsilon_0 \gamma_e + \varepsilon_0 \Delta t^2 \omega_{ep}^2 \quad (25)$$

$$c_1 = \frac{2\mu_0 \Delta t^2 \omega_{mp}^2 - 8\mu_0}{C} \quad (26)$$

$$c_2 = \frac{4\mu_0 - 2\Delta t \mu_0 \gamma_m + \mu_0 \Delta t^2 \omega_{mp}^2}{C} \quad (27)$$

$$d_0 = \frac{4 + 2\Delta t \gamma_m}{C} \quad (28)$$

$$d_1 = \frac{-8}{C} \quad (29)$$

$$d_2 = \frac{4 - 2\Delta t \gamma_m}{C} \quad (30)$$

$$C = 4\mu_0 + 2\Delta t \mu_0 \gamma_m + \mu_0 \Delta t^2 \omega_{mp}^2 \quad (31)$$

Therefore, our proposed dispersive vector meshless method is formulated based on (10), (19), (9) and (18), respectively.

### III. NUMERICAL EXAMPLE

In this section, we investigate the efficiency of our proposed method by studying a two dimensional $5mm \times 5mm$ rectangular cavity with perfectly conducting walls. As it is depicted in Fig. 1, a part of cavity is filled by plasma. We can model the plasma medium using Drude dispersion model by substitution of parameters $\omega_{ep} = \omega_p = 10^{11}$ rad/s, $\gamma_e = \gamma_m = 0$ and $\omega_{mp} = 0$ in (11) and (12).

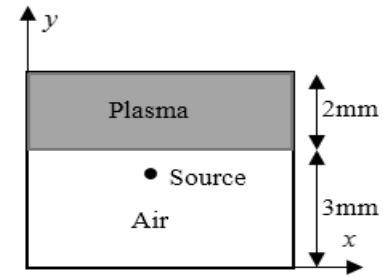

Fig. 1. Considered geometry of 2-D cavity.

A modulated Gaussian pulse is excited inside the cavity with the following form:

$$\mathbf{J}(t) = A\sin(2\pi f_0 t)\exp\left(-\left(\frac{t-t_0}{\tau}\right)^2\right)\hat{a}_y \quad (32)$$

where $f_0 = 100\ GHz, t_0 = 15\times10^{-12} s$ and $\tau = 5\times10^{-12}\ s$. Time variation of the electric field inside the cavity obtained by the proposed dispersive vector meshless method is demonstrated in Fig. 2. To investigate the numerical accuracy of the simulation, we calculated the Fourier transform of the time-domain response. As Table I shows, results are in a good agreement with the analytical resonance frequencies. Moreover, we ran simulations for the problem by dispersive scalar-based meshless method. The same nodal spacing of *0.5 mm* is considered for both methods; however, the proper shape parameters, which is attained by numerical investigation, is different for each method. We found the frequencies of two methods almost the same, but the dispersive scalar-based meshless method has a spurious mode of *25 GHz* that is far smaller than the analytical dominant frequency of the cavity. This mode is generated due to the numerical spurious charges which is calculated in the following. On the other hand, the computational time for simulations using scalar-based and vector-based dispersive meshless methods is equal to *20 s* and *72 s*, respectively. Hence, more accurate solution of vector-based method is achieved at the cost of reducing the speed of simulation.

In order to verify the divergence-free property of dispersive vector meshless, we have obtained the divergence of electric flux-density vector inside the cavity. We have located the source at the center of the cavity and computed the electric charge density, *i.e.*, $\rho$, using the following principle:

$$\rho = \nabla \cdot \mathbf{D} = \frac{\partial D_x}{\partial x} + \frac{\partial D_y}{\partial y} \quad (33)$$

The electric charge distribution is calculated using the dispersive scalar-based and vector-based meshless methods in Fig. 3. Based on the results, the electric charge density obtained through the dispersive vector meshless is more concentrated at the source location in comparison to the dispersive conventional meshless method.

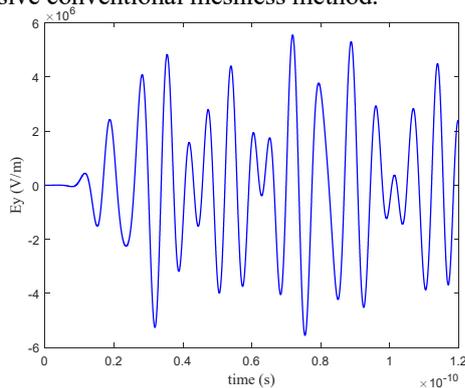

Fig. 2. $E_y$ field recorded inside the 2-D cavity solved by the dispersive vector-based meshless method.

TABLE I. RESONANCE FREQUENCY OF SIMULATION

| Analytical Results (GHz) | Dispersive Vector Meshless Method | |
|---|---|---|
| | Resonance Frequency(GHz) | Relative Error |
| 149.95 | 150 | 0.033% |
| 169.89 | 166.7 | 1.704% |

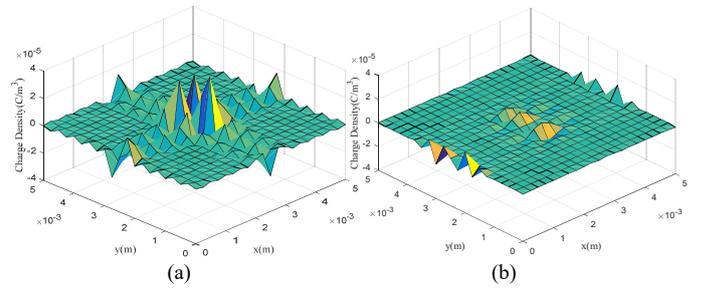

Fig. 3. The electric charge distribution obtained with (a) dispersive meshless method with scalar basis fuction and (b) dispersive meshless method with vector basis fuction at time *t* = 2 ps.

IV. CONCLUSION

In this paper, a novel dispersive vector meshless method based on ADE formulation is proposed. The advantage of dispersive vector meshless method over its scalar meshless method counterpart is its divergence-free property. The proposed method is formulated for a linear dispersive material with Drude dispersion model. Extension of the method to the other dispersive mediums can be investigated in future works.